\documentclass[11pt,english]{paper}
\usepackage{lmodern}

\usepackage[T1]{fontenc}
\usepackage[latin9]{inputenc}
\usepackage{geometry}
\usepackage{amsthm}
\usepackage{amsmath}
\usepackage{amssymb}
\usepackage{graphicx}
\usepackage{esint}

\makeatletter

\providecommand{\tabularnewline}{\\}

\numberwithin{equation}{section}
\numberwithin{figure}{section}
\numberwithin{table}{section}
\newcommand{\lyxaddress}[1]{
\par {\raggedright #1
\vspace{1.4em}
\noindent\par}
}

\makeatother

\usepackage{babel}
\begin{document}

\title{Spectral integration of linear boundary value problems }

\author{Divakar Viswanath }

\maketitle

\lyxaddress{Department of Mathematics, University of Michigan (divakar@umich.edu). }
\begin{abstract}
Spectral integration was deployed by Orszag and co-workers (1977,
1980, 1981) to obtain stable and efficient solvers for the incompressible
Navier-Stokes equation in rectangular geometries. Two methods in current
use for channel flow and plane Couette flow, namely, Kleiser-Schumann
(1980) and Kim-Moin-Moser (1977), rely on the same technique. In its
current form, the technique of spectral integration, as applied to
the Navier-Stokes equations, is dominated by rounding errors at higher
Reynolds numbers which would otherwise be within reach. In this article,
we derive a number of versions of spectral integration and explicate
their properties, with a view to extending the Kleiser-Schumann and
Kim-Moin-Moser algorithms to higher Reynolds numbers. More specifically,
we show how spectral integration matrices that are banded, but bordered
by dense rows, can be reduced to purely banded matrices. Key properties,
such as the accuracy of spectral integration even when Green's functions
are not resolved by the underlying grid, the accuracy of spectral
integration in spite of ill-conditioning of underlying linear systems,
and the accuracy of derivatives, are thoroughly explained. 
\end{abstract}

\section{Introduction}

One of the earliest methods for solving the incompressible Navier-Stokes
equation was proposed in a pioneering paper by Orszag \cite{Orszag1971}.
In that paper, Orszag tackled the problem of numerically integrating
wall-bounded shear flows using Chebyshev series expansions. The Chebyshev
polynomial is defined by $T_{n}(y)=\cos\left(n\arccos y\right)$ for
$-1\leq y\leq1$. If $u(y)=\alpha_{0}T_{0}/2+\sum_{j=1}^{\infty}\alpha_{n}T_{n}$
is the Chebyshev series of $u(y)$, we denote the Chebyshev coefficient
$\alpha_{n}$ by $\mathcal{T}_{n}(u)$. The points $y_{j}=\cos(j\pi/M)$,
$j=0,\ldots,M$, are the Chebyshev grid points. The discrete cosine
transform may be used to pass back and forth between the physical
domain function values $u(y_{j})$, $0\leq j\leq M$, and the coefficients
in the truncated Chebyshev expansion $\alpha_{0}T_{0}/2+\sum_{j=1}^{M-1}\alpha_{j}T_{j}+\alpha_{M}T_{M}/2$.

The method proposed by Orszag in \cite{Orszag1971} is certainly complete.
However, it is much too expensive. It does not appear to have been
implemented and therefore its effectiveness cannot be gaged. Nevertheless,
when numerical computations of fully turbulent solutions of shear
flows at last became possible nearly two decades later \cite{KimMoinMoser1987},
they relied on Chebyshev series expansion in the wall-normal directions
and other ideas introduced by Orszag. The two now classical methods
for computing turbulent solutions of the Navier-Stokes equation are
due to Kleiser-Schumann \cite{KleiserSchumann1980} and Kim-Moin-Moser
\cite{KimMoinMoser1987}.

The method of spectral integration was introduced by Gottlieb and
Orszag as a reformulation of the tau-equations \cite[p. 119]{GottliebOrszag1977}.
It is the thread which links the early work of Orszag \cite{Orszag1971}
with the Kleiser-Schumann and Kim-Moin-Moser methods. Below and throughout
this paper, $D$ denotes $d/dy$. The Chebyshev tau equations for
a boundary value problem such as 
\begin{equation}
\left(D^{2}-a^{2}\right)u=f(y),\, u(\pm1)=0,\label{eq:secn1-intro-bvp}
\end{equation}
are obtained by expanding the sought for solution $u$ in a truncated
Chebyshev series and equating the Chebyshev coefficients of $T_{0},\ldots,T_{M-2}$
in the expansion of $\left(D^{2}-a^{2}\right)u$ to those of $f$,
and enforcing the boundary conditions to get two more equations. As
Gottlieb and Orszag noted the tau equations are dense and not well-conditioned.
Their method of rewriting gives a tridiagonal system bordered by dense
rows corresponding to the boundary conditions. In Section 2, we derive
a variety of spectral integration methods. The Gottlieb-Orszag method
is a special case of one of them. All the methods of Section 2 work
with purely banded matrices and no bordering rows. Working with purely
banded matrices enforces a key property of spectral integration explicitly,
instead of relying on accident or the vagaries of implementation.
This key property, explained in Section 3.1, is the cancellation of
large discretization errors that arise in the intermediate stages
of spectral integration. Without this property, spectral integration
would be ineffective for the computation of turbulent solutions. The
more accurate versions of Kleiser-Schumann and Kim-Moin-Moser derived
in \cite{Viswanath2014} rely upon working with purely tridiagonal
systems.

The method of spectral integration was employed widely in the solution
of channel flow and plane Couette flow, especially in the transitional
and turbulent regimes, beginning with Orszag and coworkers \cite{OrszagKells1980,OrszagPatera1981},
Kleiser-Schumann \cite{KleiserSchumann1980}, and Kim-Moin-Moser \cite{KimMoinMoser1987}.
However, two of its key properties came to light only with the work
of Greengard \cite{Greengard1991}. The Green's function of the boundary
value problem (\ref{eq:secn1-intro-bvp}) has a scale which is $\mathcal{O}(1/a)$
and which arises from terms such as $\exp(-a|y-\tilde{y}|)$. However,
the solution may not have such a fine scale. For example, we may choose
$a=10^{6}$ and then $f(y)$ in such a way that $u(y)=\sin\pi y$
is the solution of the boundary value problem (\ref{eq:secn1-intro-bvp}).
In that situation one would need $M>2\times10^{4}$ to resolve the
Green's function at the boundaries while $M=32$ suffices to resolve
the solution $u(y)$. Greengard noted accurate solutions may be computed
when the Chebyshev grid resolves the solution even if it fails to
resolve the Green's function. 

This property appears essential for the robustness, if not the success,
of spectral integration in turbulence computations. The analogue of
$a$ in (\ref{eq:secn1-intro-bvp}) in a channel flow or plane Couette
flow computation is given by $a=\beta=\left(\gamma Re/\Delta t+\beta^{\ast}\right)^{1/2}$,
where $\beta^{\ast}$ is positive and much smaller than $Re/\Delta t$.
The Reynolds number is denoted by $Re$ and the time step by $\Delta t$;
$\gamma$ is an $\mathcal{O}\left(1\right)$ parameter. The channel
flow simulation in \cite{Viswanath2014} has $\beta\approx2\times10^{4}$
and reaches $Re=80,000$. Reaching an $Re$ that is $5$ times as
high will imply $\beta>10^{5}$ in the boundary value problems that
arise while time stepping channel flow. The solutions themselves may
not have scales as small as $\mathcal{O}\left(1/\beta\right)$, because
the smallness of $1/\beta$ is partly a reflection of time stepping
and not entirely due to the physics of the problem. Thus it is fortunate
that spectral integration can compute accurate solutions without resolving
the $1/\beta$ scale. A complete explanation of why spectral integration
has this property is given in Section 3.1. The explanation relies
on formulations in Section 2 which avoid bordering banded systems
with dense rows. 

Another property brought to light by Greengard \cite{Greengard1991}
is that condition numbers of spectral integration matrices, corresponding
to boundary value problems such as (\ref{eq:secn1-intro-bvp}), are
bounded in the limit $M\rightarrow\infty$. As noted by Rokhlin \cite{Rokhlin1983},
any integral formulation has this property because the integral operators
that are discretized are compact. In contrast, the tau equations discretize
(\ref{eq:secn1-intro-bvp}) in its differential form and therefore
suffer from ill-conditioning. In particular, their condition number
goes to $\infty$ as $M\rightarrow\infty$.

Although this is a useful property, it is by itself inadequate to
understand the robustness of spectral integration as applied to the
Navier-Stokes equations in the turbulent regime. If $a=10^{6}$ in
the boundary value problem (\ref{eq:secn1-intro-bvp}), for example,
the condition number of the spectral integration matrices is of the
order $\mathcal{O}\left(10^{12}\right)$. For $4$-th order problems
such as $\left(D^{2}-a^{2}\right)\left(D^{2}-b^{2}\right)u=f$, the
condition number appear to get as large as $\mathcal{O}\left(a^{2}b^{2}\right)$
in the limit $M\rightarrow\infty$. The fact that the condition numbers
do not diverge as $M\rightarrow\infty$ offers little comfort when
condition numbers become so large. Such large condition numbers suggest
a severe loss of accuracy. Yet even in the presence of such ill-conditioning,
spectral integration is able to compute solutions $u(y)$ with close
to machine precision. In Section 3.2, we explain how spectral integration
comes to have such remarkable accuracy. 

Greengard's form of spectral integration as applied to (\ref{eq:secn1-intro-bvp})
expands $D^{2}u$ instead of $u$ in a Chebyshev series. Thus it may
be suspected that Greengard's form produces more accurate derivatives.
In fact, that is partially but not entirely true. Building upon a
discovery of Muite \cite{Muite2010}, we discuss why all forms of
spectral integration produce derivatives with similar accuracy when
the parameter $a$ of (\ref{eq:secn1-intro-bvp}) is not too large.
The key point, as noted by Muite, is that one must not pass into physical
space when calculating derivatives. However, when $a$ is large, as
in turbulence simulations with high Reynolds number, there is a significant
advantage to avoiding numerical differentiation altogether.

The more robust versions of Kleiser-Schumann and Kim-Moin-Moser derived
in \cite{Viswanath2014} are specially adapted to either method and
combine several of the forms of spectral integration derived in Section
2. They altogether avoid numerical differentiation in the wall-normal
direction. The ease with which the different forms can be combined
follows from the way we enforce boundary conditions. Instead of using
ad-hoc bordered rows for boundary conditions, we derive the boundary
conditions from integral conditions that are intrinsic to spectral
integration formulations. Both the method of Gottlieb and Orszag \cite[p. 119]{GottliebOrszag1977}
and that of Greengard \cite{Greengard1991} are seen to be special
cases. 

Certain comments made by Orszag and coworkers \cite{OrszagKells1980,OrszagPatera1981}
show awareness of the connection of the method of Gottlieb and Orszag
\cite[p. 119]{GottliebOrszag1977} to integral formulations. A more
explicit connection was made by Muite \cite{Muite2010} and by Charalambides
and Waleffe \cite{CharalambidesWaleffe2008A,CharalambidesWaleffe2008B}.
The latter authors study spurious eigenvalues of spectral discretization
matrices using the theory of stable polynomials. They consider Jacobi
and Gegenbauer polynomials in addition to Chebyshev. A connection
of the method of Gottlieb and Orszag to spectral integration is made
at a more fundamental level in Section 2.

A number of numerical examples are included in Section 4. These examples
are illustrative but do not indicate the potential of the methods
derived here. When problem sizes are small, it is difficult to argue
for the value of carefully derived and highly accurate methods. Therefore
Section 4 discusses the Kleiser-Schumann algorithm in the context
of properties of spectral integration derived in Section 2. The main
application is to the integration of the Navier-Stokes equation in
rectangular geometry and in the turbulent regime. In a companion paper
\cite{Viswanath2014}, we present a computation of turbulent channel
flow that reaches $Re_{\tau}=2380$ using $10^{9}$ grid points and
$M=1024$. Only $10$ nodes of a small cluster are used in this computation.
The computation uses an algorithm, derived using the techniques of
the present paper, which avoids numerical differentiation in the wall-normal
direction entirely. This computation may be compared with those of
Hoyas and Jiménez \cite{HoyasJimenez2006,HoyasJimenez2008} who reached
$Re_{\tau}=2003$. As explained in \cite{Viswanath2014}, current
methods have fundamental limitations that limit them to $M\approx1000$.
The methods derived in \cite{Viswanath2014} using techniques of the
present paper appear capable of handling problems with $M=10^{4}$,
which would take us into problem sizes beyond the reach of modern
computers. Turbulence computations such as the ones in \cite{HoyasJimenez2006,HoyasJimenez2008,Viswanath2014}
are typically run for several months. 

Beyond the application to channel flow and plane Couette flow, the
methods derived here are likely to be of use in other problems in
fluid mechanics such as Rayleigh-Bénard convection and Kolmogorov
flow.

\section{Varieties of spectral integration}

In this section, the interval of the boundary value problem is taken
to be $-1\leq y\leq1$ and the solution $u$ or one of its derivatives
is expanded as follows:
\begin{equation}
\frac{\alpha_{0}}{2}+\alpha_{1}T_{1}+\cdots+\alpha_{M-1}T_{M-1}+0.T_{M}.\label{eq:u-cheb-expansion}
\end{equation}
The Chebyshev points $y_{j}=\cos(j\pi/M)$ with $j=0\ldots M$ are
$M+1$ in number including the endpoints $\pm1$, but the last coefficient
in the Chebyshev series is suppressed for convenience as indicated
in (\ref{eq:u-cheb-expansion}). 

In each of the methods of this section including the factored form
of spectral integration, the way the homogeneous solutions are computed
may appear roundabout. As explained in Section 3.1, such a roundabout
calculation is essential for ensuring accuracy. 

The method of Gottlieb and Orszag \cite[p. 119]{GottliebOrszag1977}
is a special case of the method in Section 2.1. The method of Greengard
\cite{Greengard1991} is a special case of the method in Section 2.3.
Apart from being more general, our formulation works with purely banded
systems instead of banded systems with rows. If the boundary conditions
are not suitably scaled, Gaussian elimination with partial pivoting
can turn a banded system with a few dense rows into a dense matrix.
The use of QR factorization to solve such systems has been proposed
\cite{OlverTownsend2013}. QR with its use of square roots would be
much too slow for applications such as computations of turbulence
in channel flow. All these issues are swept aside by the methods given
here, all of which use only purely banded systems. Another advantage
is that our formulations lead to an explanation of why spectral integration
schemes are accurate if they resolve the solution but not the Green's
function. This explanation is given in Section 3.1.

\subsection{First and second order spectral integration}

Suppose a first order boundary value problem is given in the form
\begin{equation}
(D-a)u=f.\label{eq:1st-order-template-problem}
\end{equation}
The exact boundary condition is unimportant for much of the method.
To begin with we assume the integral condition $\mathcal{T}_{0}(u)=0$
or equivalently $\alpha_{0}=0$. Suppose that the Chebyshev coefficients
of $f$ are given by $f_{j}=\mathcal{T}_{j}(f)$. 

The indefinite integral of (\ref{eq:1st-order-template-problem})
gives 
\begin{equation}
u-a\int u+A=\int f,\label{eq:first-order-template-integrated}
\end{equation}
where $A$ is an undetermined constant. The integral $\int T_{n}\, dy$
is $T_{n+1}/2(n+1)-T_{n-1}/2(n-1)$ if $n>1$, $T_{2}/4$ if $n=1$,
and $T_{1}$ if $n=0$. Therefore, the coefficient of $T_{n}$ on
the right hand side of (\ref{eq:first-order-template-integrated})
is 
\[
\mathcal{T}_{n}\left(\int f\right)=\frac{f_{n-1}-f_{n+1}}{2n}\quad\text{for}\quad n=1,2,\ldots,M-1.
\]
The $n=M-1$ case assumes $f_{M}=0$. Similarly, the coefficient of
$T_{n}$ in the expansion of the left hand side of (\ref{eq:first-order-template-integrated})
is 
\[
\mathcal{T}_{n}\left(u-a\int u+A\right)=\begin{cases}
\alpha_{n}-a\left(\frac{\alpha_{n-1}-\alpha_{n+1}}{2n}\right) & 2\leq n\leq M-2\\
\alpha_{n}-a\left(\frac{-\alpha_{n+1}}{2n}\right) & n=1\\
\alpha_{n}-a\left(\frac{\alpha_{n-1}}{2n}\right) & n=M-1
\end{cases}
\]
for $n=1,\ldots,M-1$. The coefficients with $n=1$ and $n=M-1$ are
obtained from the more general expression $\alpha_{n}-a(\alpha_{n-1}-\alpha_{n+1})/2n$
by setting $\alpha_{0}=0$ and $\alpha_{M}=0$, respectively. Equating
coefficients for $n=1,\ldots,M-1$ we have $M-1$ equations for the
$M-1$ unknowns $\alpha_{1},\ldots,\alpha_{M-1}$. If we did not set
$\alpha_{M}=0$, another equation with $n=M$ may be used, but that
equation has a different form from the equations for $1\leq n\leq M-1$.
Setting $\alpha_{M}=0$ saves us from a little inconvenience. This
tridiagonal linear system is solved to compute a particular solution
$u^{p}$ satisfying $(D-a)u^{p}=f$ and the integral condition $\mathcal{T}_{0}(u^{p})=0$.

A homogeneous solution $\bar{u}^{1}$ satisfying $(D-a)\bar{u}^{1}=0$
and $\mathcal{T}_{0}(\bar{u}^{1})=1$ is found as follows. We set
$\bar{u}=1/2+u^{\ast}$ so that $u^{\ast}$ satisfies $\mathcal{T}_{0}(u^{\ast})=0$
and $(D-a)u^{\ast}=a/2$. Thus $u^{\ast}$ is the particular solution
of (\ref{eq:1st-order-template-problem}) satisfying $\mathcal{T}_{0}(u^{\ast})=0$
if $f\equiv-a/2$ and it may be found using the method described for
computing particular solutions. The same linear tridiagonal system
is solved for computing $u^{\ast}$ and $u^{p}$ but with different
right hand sides.

The solution $u$ is expressed as $u^{p}+C\bar{u}^{1}$ and the constant
$C$ is found using the boundary condition on $u$.

Now we consider the second order problem 
\begin{equation}
\left(D^{2}+bD+c\right)u=f.\label{eq:second-order-template}
\end{equation}
Integrating twice assuming $b,c$ to be constant, we have
\[
u+b\int u+c\int\int u+A+By=\int\int f.
\]
To find a particular solution, we assume the integral conditions $\mathcal{T}_{0}(u)=\mathcal{T}_{1}(u)=0$
or equivalently $\alpha_{0}=\alpha_{1}=0$. By standard formulas for
$\int T_{n}$ and $\int\int T_{n}$, the coefficient of $T_{n}$ of
the right hand side is
\[
\mathcal{T}_{n}\left(\int\int f\right)=\frac{f_{n-2}}{4n(n-1)}-\frac{f_{n}}{2(n^{2}-1)}+\frac{f_{n+2}}{4n(n+1)}
\]
for $n=2,\ldots,M-1$ (here $f_{M}=f_{M+1}=0$ is assumed). The coefficient
of the left hand side is
\begin{multline*}
\mathcal{T}_{n}\left(u+b\int u+c\int\int u+A+By\right)=\alpha_{n-2}\left(\frac{c}{4n(n-1)}\right)+\alpha_{n-1}\left(\frac{b}{2n}\right)\\
+\alpha_{n}\left(1-\frac{c}{2(n^{2}-1)}\right)+\alpha_{n+1}\left(\frac{-b}{2n}\right)+\alpha_{n+2}\left(\frac{c}{4n(n+1)}\right)
\end{multline*}
for $n=2,\ldots,M-1$. The validity of the equations for $n=2,3$
relies on the boundary conditions $\alpha_{1}=\alpha_{2}=0$. The
equations for $n=M-2,M-1$ assume $\alpha_{M}=\alpha_{M+1}=0$. The
coefficients for $n=2,\ldots,M-1$ on the left and right hand sides
are equated to solve for the unknowns $\alpha_{2},\ldots,\alpha_{M-1}$.
The particular solution $u^{p}$ obtained in this manner satisfies
$(D^{2}+bD+c)u^{p}=f$ and the integral conditions $\mathcal{T}_{0}(u^{p})=\mathcal{T}_{1}(u^{p})=0$.

The first homogeneous solution satisfies $(D^{2}+bD+c)\bar{u}^{1}=0$
and the integral conditions $\mathcal{T}_{0}(\bar{u}^{1})=1,\mathcal{T}_{1}(\bar{u}^{1})=0$.
To find it, we set $\bar{u}^{1}=1/2+u^{\ast}$. Then $u^{\ast}$ satisfies
the inhomogeneous equation (\ref{eq:second-order-template}) with
$f\equiv-c/2$ and the integral conditions $\mathcal{T}_{0}(u^{\ast})=\mathcal{T}_{1}(u^{\ast})=0$.
The solution $u^{\ast}$ is computed using the same pentadiagonal
system used for $u^{p}$ but with a different right hand side.

The second homogeneous solution $\bar{u}^{2}$ satisfies the integral
conditions $\mathcal{T}_{0}(\bar{u}^{2})=0,\mathcal{T}_{1}(\bar{u}^{2})=1$.
If we set $\bar{u}^{2}=T_{1}+u^{\ast}$, $u^{\ast}$ satisfies the
inhomogeneous equation with $f\equiv-(b+cT_{1})$. It is found by
solving the same pentadiagonal system.

The solution $u$ of (\ref{eq:second-order-template}) is expressed
as $u^{p}+C\bar{u}^{1}+D\bar{u}^{2}$ and the constants $C$ and $D$
are determined using the boundary conditions on $u$.

\subsection{Spectral integration of $r$-th order}

Define the operator $L$ as $Lu=u^{(r)}+a_{r-1}u^{(r-1)}+\cdots+a_{1}u^{(1)}+a_{0}u$.
Consider the inhomogeneous equation $Lu=f$. A particular solution
satisfying the integral conditions 
\[
\mathcal{T}_{0}(u)=\cdots=\mathcal{T}_{r-1}(u)=0
\]
or equivalently $\alpha_{0}=\cdots=\alpha_{r-1}=0$ may be found as
follows. Assuming constant coefficients, the inhomogeneous equation
is written in an integral form as 
\[
u+a_{r-1}\int u+\cdots+a_{0}\int^{r}u+\sum_{j=0}^{r-1}A_{j}y^{j}=\int^{r}f.
\]
Using formulas for $\int^{j}T_{n}$ we may express the coefficients
of $T_{r},\ldots,T_{M-1}$ in terms of $\alpha_{r},\ldots,\alpha_{M-1}$.
These coefficients are equated to the coefficients of $\int^{r}f$
and solved for $\alpha_{r},\ldots,\alpha_{M-1}$ to find the particular
solution $u^{p}$. The linear system has $2r+1$ diagonals.

To find the $j$-th homogeneous solution for $j=1,\ldots,r$, we first
set $\bar{u}^{j}=T_{j-1}+u^{\ast}$. The $j$-th homogeneous solution
satisfies the conditions $\mathcal{T}_{k}(\bar{u}^{j})=0$, if $0\leq k\leq r-1$
and $k\neq j-1$, and $\mathcal{T}_{j-1}(\bar{u}^{j})=1$. The function
$u^{\ast}$ satisfies $Lu^{\ast}=-LT_{j}$ and the first $r$ coefficients
in its Chebyshev series are zero. It can be found in the same manner
as the particular solution.

The solution of the linear boundary value problem is expressed as
$u^{p}+\sum_{j=1}^{r}C_{j}\bar{u}^{j}$. The boundary conditions satisfied
by $u$ are used to determine the constants $C_{j}$.

Formulas for $\int^{j}T_{n}$ can be derived but get complicated.
The $2r+1$ diagonal system can be difficult to set up correctly in
programs. For the difficulties that arise for $r=4$, see \cite{Muite2010}.
Spectral integration of $r$-th order will, however, prove quite useful
in the discussion of cancellation errors in Section 3.1.

\subsection{Greengard form of spectral integration}

The formulation of the Greengard form of spectral integration given
here makes it clear that the matrix systems solved have $2r+1$ diagonals
if the problem is of order $r$. We assume $L$ to be the operator
defined in Section 2.2 above.

The Greengard form begins by assuming a Chebyshev series for $u^{(r)}$.
A similar method was proposed earlier by Zebib \cite{Zebib1984}.
We first find a particular solution of $Lu=f$ subject to the integral
conditions 
\[
\mathcal{T}_{0}(u)=\mathcal{T}_{0}(u^{(1)})=\cdots=\mathcal{T}_{0}(u^{(r-1)})=0.
\]
The integral conditions are different this time. The integral conditions
given earlier ensure that if $u$ satisfies the integral conditions
and we know the Chebyshev series of $u$, then we can produce the
Chebyshev series of $\int_{s}u$ for $s=0,\ldots r-1$. The integral
conditions here ensure that if we know the Chebyshev series of $u^{(r)}$,
then we can produce the Chebyshev series of $u^{(s)}$ for $s=0,\ldots,r-1$
without ambiguity. When these conditions are used, the Chebyshev series
of $u^{(s)}$ determines the Chebyshev series of $u^{(s-1)}$ for
$s=r,\ldots,1$ with no ambiguity. Normally, there is an undetermined
constant of integration when the series of $u^{(s)}$ is integrated.
But here the constant disappears because the mean mode of $u^{(s-1)}$
is specified to be zero. Thus the Chebyshev series of $Lu$ is determined
unambiguously by the Chebyshev series of $u^{(r)}$. Coefficients
in the Chebyshev series of $Lu$ and $f$ are equated and solved for
$\mathcal{T}_{j}\left(u^{(r)}\right)$ for $j=0,\ldots,M$.

To find homogeneous solutions $\bar{u}$, we expand $\bar{u}^{(r)}$
in a Chebyshev series and take the integral conditions to be such
that exactly one of $\mathcal{T}_{0}(\bar{u}),\ldots,\mathcal{T}_{0}(\bar{u}^{(r-1)})$
is one and the others are all zero. It is harder to find homogeneous
solutions here than in the $r$-th order spectral integration method
described in Section 2.2. One has to find polynomials $p_{k}$ of
degree $k$ for $k=0,1,\ldots,r-1$ such that $\mathcal{T}_{0}(p_{k}^{(k)})=1$
but $\mathcal{T}_{0}(p_{k}^{(d)})=0$ for $d=0,\ldots,k-1$. This
form of spectral integration generalizes more easily to linear differential
equations with polynomial coefficients \cite{CoutsiasHT1996}.

\subsection{Factored form of spectral integration}

A linear operator $L$ with constant and real coefficients can be
factorized as 
\[
L=(D-a_{1})\ldots(D-a_{m})(D^{2}+b_{1}D+c_{1})\ldots(D^{2}+b_{n}D+c_{n})
\]
where the coefficients are all real. We assume $m+n\geq2$ and derive
a method for solving $Lu=f$ subject to boundary conditions that exploits
this factorization of $L$. This method relies on spectral integration
of orders one and two described in Section 2.1. The presentation of
the method may appear more complicated but its implementation is much
simpler than the methods of Sections 2.2 and 2.3. 

A particular solution is found by solving the following equations
subject to integral conditions on their solutions:
\begin{gather*}
(D-a_{1})u_{m+2n-1}^{p}=f\quad\quad\mathcal{T}_{0}(u_{m+2n-1}^{p})=0\\
(D-a_{2})u_{m+2n-2}^{p}=u_{m+2n-1}^{p}\quad\quad\mathcal{T}_{0}(u_{m+2n-2}^{p})=0\\
\vdots\\
(D-a_{m})u_{2n}^{p}=u_{2n+1}^{p}\quad\quad\mathcal{T}_{0}(u_{2n}^{p})=0\\
(D^{2}+b_{1}D+c_{1})u_{2n-2}^{p}=u_{2n}^{p}\quad\quad\mathcal{T}_{0}(u_{2n-2}^{p})=\mathcal{T}_{1}(u_{2n-2}^{p})=0\\
\vdots\\
(D^{2}+b_{n}D+c_{n})u_{0}^{p}=u_{2}^{p}\quad\quad\mathcal{T}_{0}(u_{0}^{p})=\mathcal{T}_{1}(u_{0}^{p})=0.
\end{gather*}
This list of equations is solved from first to last. Each equation
is solved using one of the two methods described in Section 2.1. The
subscripts on $u$, as in $u_{2n}^{p}$ , indicate the number of ``derivatives''
in the function relative to $u_{0}^{p}$ which satisfies $Lu_{0}^{p}=f$
and is therefore a particular solution.

If $m\geq1$, the homogeneous solution $\bar{u}^{h}$ with $h=1$
is found as follows. To begin with we solve the homogeneous problem
\[
(D-a_{1})\bar{u}_{n+2m-1}^{h}=0\quad\quad\mathcal{T}_{0}(\bar{u}_{n+2m-1}^{h})=1
\]
as described in Section 2.1. Thereafter, the inhomogeneous problems
\begin{equation}
(D-a_{j})\bar{u}_{n+2m-j}^{h}=\bar{u}_{n+2m-j+1}^{h}\quad\quad\mathcal{T}_{0}(\bar{u}_{n+2m-j}^{h})=0\label{eq:intermed-aj}
\end{equation}
are solved in the order $j=2,\ldots,m$ followed by the solution of
\begin{equation}
(D^{2}+b_{k}D+c_{k})\bar{u}_{2n-2k}^{h}=\bar{u}_{2n-2k+2}^{h}\quad\quad\mathcal{T}_{0}(\bar{u}_{2n-2k}^{h})=\mathcal{T}_{1}(\bar{u}_{2n-2k}^{h})=0\label{eq:intermed-bkck}
\end{equation}
in the order $k=1,\ldots,n$. The last solution to be found is $\bar{u}^{h}=\bar{u}_{0}^{h}$
and it satisfies $L\bar{u}^{h}=0$. The inhomogeneous equations (\ref{eq:intermed-aj})
and (\ref{eq:intermed-bkck}) are solved as described in Section 2.1.

More generally, the homogeneous solution $\bar{u}^{h}$ with $1\leq h\leq m$
is solved beginning with the homogeneous problem 
\[
(D^{2}-a_{h})\bar{u}_{m+2n-h}^{h}=0\quad\quad\mathcal{T}_{0}(\bar{u}_{m+2n-h}^{h})=1
\]
followed by the solution of (\ref{eq:intermed-aj}) with $j=h+1,\ldots,m$
and (\ref{eq:intermed-bkck}) with $k=1,\ldots,n$. As before, $\bar{u}_{0}^{h}$
is the last solution to be found and $\bar{u}^{h}=\bar{u}_{0}^{h}$.

If $h=m+2i-1$ with $1\leq i\leq n$, the homogeneous problem solved
at the beginning is 
\[
(D^{2}+b_{i}D+c_{i})\bar{u}_{2n-2i}^{h}=0\quad\quad\mathcal{T}_{0}(\bar{u}_{2n-2i}^{h})=1,\,\mathcal{T}_{1}(\bar{u}_{2n-2i}^{h})=0.
\]
This is followed by the solution of (\ref{eq:intermed-bkck}) with
$k=i+1,\ldots,n$. As before, $\bar{u}_{0}^{h}$ is the last solution
to be found and $\bar{u}^{h}=\bar{u}_{0}^{h}$. On the other hand,
if $h=m+2i$ with $1\leq i\leq n$, the homogeneous problem solved
at the beginning is 
\[
(D^{2}+b_{i}D+c_{i})\bar{u}_{2n-2i}^{h}=0\quad\quad\mathcal{T}_{0}(\bar{u}_{2n-2i}^{h})=0,\,\mathcal{T}_{1}(\bar{u}_{2n-2i}^{h})=1.
\]
This is followed by the solution of (\ref{eq:intermed-bkck}) with
$k=i+1,\ldots,n$. As before, $\bar{u}_{0}^{h}$ is the last solution
to be found and $\bar{u}^{h}=\bar{u}_{0}^{h}$.

By using the methods of Section 2.1 repeatedly, we end up with a particular
solution $u^{p}$ and homogeneous solutions $\bar{u}^{1},\ldots,\bar{u}^{m+2n}$.
The solution of the boundary value problem $Lu=f$ is expressed as
\[
u=u^{p}+\sum_{j=1}^{m+2n}C_{j}\bar{u}^{j}.
\]
The constants $C_{j}$ are found to fit the boundary conditions on
$u$. 

There are two ways to find $C_{j}$. In the first method, the particular
solution $u^{p}$ and the homogeneous solutions $\bar{u}^{h}$ are
obtained in physical space as numerical values at the $M+1$ points
on the Chebyshev grid. Boundary conditions such as $u(1)=A$ or $u''(-1)=B$
are expressed using a linear combinations of function values at the
grid point. A boundary condition such as $u(1)=A$ simply specifies
the function value at a single grid point. A boundary condition such
as $u''(-1)$ is interpreted as specifying that a certain linear combination
of function values, the linear combination being determined by a single
row of a spectral differentiation matrix, must have a specified value.
This is the easier method for implementation and the one we have implemented.
If the number $M$ is not too large, this method will be adequate.
If $M$ is very large, then errors will creep in through the boundaries. 

The second technique uses the intermediate objects created when the
particular solution and the homogeneous solutions are found. We illustrate
the technique using an example. Suppose the boundary value problem
is 
\[
(D-a_{1})(D-a_{2})(D-a_{3})(D-a_{4})u=f
\]
subject to $u(\pm1)=u'(\pm1)=0$. If $u=u^{p}+\sum_{j=1}^{4}C_{j}\bar{u}^{j}$,
the conditions on $u(\pm1)$ give two equations for the $C_{j}$ after
evaluation at $\pm1$. We may rewrite the other boundary conditions
as $(D-a_{4})u=0$ at $\pm1$. If we now note that 
\[
(D-a_{4})u=u_{1}^{p}+\sum_{j=1}^{3}C_{j}\bar{u}_{1}^{j}
\]
we get two more equations for the $C_{j}$ by evaluating at $\pm1$.
In the sum above, the $j=4$ term does not appear. That is because
the homogeneous solution $\bar{u}^{4}$ satisfies $(D-a_{4})\bar{u}^{4}=0$.

In light of the discussion given here, a part of the methods in \cite{KimMoinMoser1987,KleiserSchumann1980,OrszagPatera1981}
may be viewed as special cases of the factored form of spectral integration.
The treatment given here suggests the more powerful versions of Kleiser-Schumann
and Kim-Moin-Moser that are derived in \cite{Viswanath2014}. The
forms of spectral integration in Sections 2.1, 2.2, and 2.3 can be
combined seamlessly with spectral integration in its factored from.
This flexibility proves to be useful in deriving versions of Kleiser-Schumann
and Kim-Moin-Moser that are resistant to rounding errors.

\subsection{Spectral integration with piecewise Chebyshev grid}

To generalize spectral integration to piecewise Chebyshev grids, we
consider the operator $(D-a)(D-b)$ over the interval $-1\leq y\leq1$
and the boundary value problem corresponding to $(D-a)(D-b)u=f$.
As earlier in this section, the boundary conditions enter only at
the end and much of the method is independent of the specific form
of the boundary conditions. The generalization to operators of the
form $(D-a_{1})\ldots(D-a_{m})$ will be obvious.

Let $[-1,1]=\mathcal{I}_{1}\cup\ldots\cup\mathcal{I}_{n}$, where
$\mathcal{I}_{i}$ are intervals with disjoint interiors and with
the right end point of $\mathcal{I}_{i}$ equal to the left end point
of $\mathcal{I}_{i+1}$. Thus $\mathcal{I}_{1},\ldots,\mathcal{I}_{n}$
are disjoint intervals arranged in order. Let $w_{i}$ denote the
width of the interval $\mathcal{I}_{i}$.

We use a linear change of variables $\mathcal{I}_{i}\rightarrow[-1,1]$
and rewrite the given differential equation as 
\begin{equation}
\left(D-\frac{aw_{i}}{2}\right)\left(D-\frac{bw_{i}}{2}\right)u=\frac{w_{i}^{2}}{2}\, f\label{eq:DaDb-rescaled}
\end{equation}
after the change of variables. In (\ref{eq:DaDb-rescaled}) it is
assumed that $u$ and $f$ have been shifted from $\mathcal{I}_{i}$
to $[-1,1]$ although that is not indicated explicitly by the notation. 

We define $u_{i}$ as 
\begin{equation}
u_{i}=u^{p}+\alpha_{\mathcal{I}_{i}}\bar{u}^{1}+\beta_{\mathcal{I}_{i}}\bar{u}^{2},\label{eq:ui-defn}
\end{equation}
 where $u^{p}$ is the particular solution and $\bar{u}^{1},\,\bar{u}^{2}$
are the homogeneous solutions of (\ref{eq:DaDb-rescaled}), computed
as described in Section 3. 

For $i=1,\ldots,n$, the coefficients $\alpha_{\mathcal{I}_{i}}$
and $\beta_{\mathcal{I}_{i}}$ comprise $2n$ unknown variables in
total. We will solve for these unknowns using the two boundary conditions
and continuity conditions between intervals. The boundary conditions
give two equations such as 
\begin{align*}
u_{1}(-1) & =\text{left value}\\
u_{n}(1) & =\text{right value}.
\end{align*}
For $i=1,\ldots,n-1$, the continuity conditions are 
\begin{align*}
u_{i}(1) & =u_{i+1}(-1)\\
\frac{\left((D-w_{i}b/2)u_{i}\right)(1)}{w_{i}} & =\frac{\left((D-w_{i+1}b/2)u_{i+1}\right)(-1)}{w_{i+1}}.
\end{align*}
The second continuity condition requires the derivatives to be continuous
while accounting for the shifting and scaling of intervals of width
$w_{i}$ and $w_{i+1}$ to $[-1,1]$. The function $(D-w_{i}b/2)u_{i}$
is available through the intermediate quantities generated by the
method of Section 3. In particular, we have, 
\[
(D-w_{i}b/2)u^{p}=u_{1}^{p},\quad(D-w_{i}b/2)\bar{u}^{1}=\bar{u}_{1}^{1},\quad(D-w_{i}b/2)\bar{u}^{2}=0
\]
in interval $\mathcal{I}_{i}$. Once we solve for $\alpha_{\mathcal{I}_{i}}$
and $\beta_{\mathcal{I}_{i}}$ for $i=1,\ldots,n$, we may use (\ref{eq:ui-defn})
to form $u_{i}$. The solution $u$ is obtained by shifting the $u_{i}$
from $[-1,1]$ back to $\mathcal{I}_{i}$.

It is important to note that the system of equations for finding $\alpha_{\mathcal{I}_{i}}$
and $\beta_{\mathcal{I}_{i}}$ is banded. The use of banded matrices
is an improvement of the gluing procedure of Greengard and Rokhlin
\cite{GreengardRokhlin1991}, although its applicability is more limited.

\section{Properties of spectral integration}

\subsection{Cancellation of intermediate errors}

The solution of the linear boundary value problem $(D^{2}-a^{2})u=-(\pi^{2}+a^{2})\sin\pi y$
with boundary conditions $u(\pm1)=0$ is $u=\sin\pi y$. The solution
can be represented with machine precision on a Chebyshev grid that
uses slightly more than $20$ points. If $a=10^{6}$ it will take
a Chebyshev grid with $M>2\times10^{4}$ points to resolve the Green's
function at the boundaries. Spectral integration can solve this boundary
value problem using a Chebyshev grid with $20$ or $30$ points even
if $a=10^{6}$. In this section, we explain the rather roundabout
manner in which spectral integration comes to acquire this property.

If $x$ is the solution of the matrix system $Ax=b_{1}+b_{2}$ and
$x_{1},\: x_{2}$ are the solutions of $Ax_{1}=b_{1},\: Ax_{2}=b_{2}$,
then $x=x_{1}+x_{2}$ if $A$ is nonsingular. In machine arithmetic
and in the presence of rounding errors, this linear superposition
property will be true only approximately. This section deals with
discretization errors and not rounding errors. Therefore we will assume
this linear superposition property.

Suppose that $Lu=f$ is the given equation. With given boundary conditions
on $u$, this equation is assumed to have a solution that is well-resolved
using $M+1$ Chebyshev points. In all forms of spectral integration,
a particular solution satisfying $Lu^{p}=f$ is found using some other
global conditions on $u^{p}$. Typically it will take many more points
than $M$ to resolve $u^{p}$. Thus the computed $u^{p}$ will be
inaccurate. However, the approximation to $u$ obtained by combining
$u^{p}$ with the homogeneous solutions will retain its accuracy for
reasons we will now explain.

The explanation takes its simplest form for $r$-th order spectral
integration described in Section 2.2 and it is with that method that
we begin. Suppose $L$ is a linear differential operator with constant
coefficients and order $r$ as in Section 2.2. We begin by denoting
the \emph{computed} solution of 
\[
Lu=LT_{j}
\]
with integral conditions $\mathcal{T}_{0}(u)=\cdots=\mathcal{T}_{r-1}(u)=0$
by $U_{j}$ for $j=0,\ldots,r-1$. Thus the Chebyshev series of $U_{j}$
is obtained by solving a banded system with $2r+1$ diagonals and
a right hand side that corresponds to the Chebyshev series of $LT_{j}$
integrated $r$ times. The $U_{j}$ will be typically quite inaccurate.
We will show that the $U_{j}$ occur in the particular solution and
the homogeneous solutions in such a way that they cancel when an approximation
to the solution of $Lu=f$ with the given boundary conditions is computed.

Let $u_{E}$ be the solution of $Lu_{E}=f$ which satisfies the given
boundary conditions and is accurate to machine precision with a Chebyshev
series of $M$ terms.We rewrite $u_{E}$ as 
\[
u_{E}=\frac{\alpha_{0}}{2}+\alpha_{1}T_{1}+\cdots+\alpha_{r-1}T_{r-1}+u_{R}
\]
where $\mathcal{T}_{j}(u_{R})=0$ for $j<r$. We may rewrite $f$
as 
\[
f=\frac{\alpha_{0}}{2}LT_{0}+\cdots+\alpha_{r-1}LT_{r-1}+f_{R}
\]
where $f_{R}=Lu_{R}$.

The particular solution of $Lu=f_{R}$ which is computed by $r$-th
order spectral integration is $u^{p}=u_{R}$. This is because $u_{R}$
satisfies the integral boundary conditions, the first $r$ of its
Chebyshev coefficients being zero, as well as $Lu=f_{R}$ and can
be represented to machine precision using a Chebyshev series of $M$
terms. By linear superposition, the particular solution of $Lu=f$
satisfying integral boundary conditions that is computed is given
by 
\begin{equation}
u^{p}=\frac{\alpha_{0}}{2}U_{0}+\alpha_{1}U_{1}+\cdots+\alpha_{r-1}U_{r-1}+u_{R}.\label{eq:computed-up}
\end{equation}

Homogeneous solutions of $Lu=0$ are computed such that $\mathcal{T}_{j}(u)=1$
but with the other $r-1$ Chebyshev coefficients among the first $r$
coefficients being zero. This homogeneous solution is represented
as $u=T_{j}+u^{\ast}$ and $u^{\ast}$ is computed as the particular
solution of $Lu^{\ast}=-LT_{j}$, whose first $r$ coefficients are
zero. Therefore the computed homogeneous solutions are 
\begin{equation}
\bar{u}^{1}=1/2-U_{0}/2,\quad\bar{u}^{2}=T_{1}-U_{1},\quad\ldots,\quad\bar{u}^{r}=T_{r-1}-U_{r-1}.\label{eq:computed-uh}
\end{equation}
By observing (\ref{eq:computed-up}) and (\ref{eq:computed-uh}),
we recognize that 
\[
u_{E}=u^{p}+\frac{\alpha_{0}}{2}\bar{u}^{1}+\alpha_{1}\bar{u}^{2}+\cdots+\alpha_{r-1}\bar{u}^{r}.
\]
In this linear combination of the particular solution with the homogeneous
solutions, the coefficients are such that the inaccurate $U_{j}$
cancel exactly and the solution $u_{E}$ satisfies the given boundary
conditions. If the equations that are solved to determine the linear
combination of homogeneous solutions with the particular solution
are reasonably well-conditioned, which we may expect because these
are typically very small linear systems, the computed solution will
produce $u_{E}$ very accurately.

The explanations for the factored form of spectral integration and
the Zebib-Greengard version are more complicated. We will give the
explanation for the problem $(D-a)(D-b)u=f$. The given boundary conditions
are assumed to be $u(\pm1)=0$. We assume as before that $u_{E}$
is the approximate solution whose Chebyshev series has $M$ terms
and which is accurate to machine precision. 

Suppose $U_{2}$ is the solution of $(D-b)u=1$ satisfying $\mathcal{T}_{0}(U_{2})=0$
computed as explained in Section 2.1 using a Chebyshev series with
$M$ terms. Similarly, let $U_{1}'$ be the computed solution of $(D-a)u=1$
satisfying $\mathcal{T}_{0}(U_{1}')=0$, and let $U_{1}$ be the particular
solution of $(D-b)u=U_{1}'$ satisfying $\mathcal{T}_{0}(U_{1})=0$
and computed using a Chebyshev series with $M$ terms only. For reasons
given above, $U_{1}$ and $U_{2}$ are typically very inaccurate. 

As before, we will split $u_{E}$ but the split is more complicated
this time. We write 
\[
u_{E}=\frac{\alpha_{0}}{2}+(\alpha_{1}-\gamma)T_{1}+u_{R},
\]
where $\gamma$ is chosen such that 
\[
u_{R}=\gamma T_{1}+\alpha_{2}T_{2}+\cdots+\alpha_{M-1}T_{M-1}
\]
satisfies $(D-b)u_{R}=0.$ By applying $(D-a)(D-b)$ to $u_{E}$,
$f$ can be split as 
\[
f=\frac{ab\alpha_{0}}{2}-(a+b)(\alpha_{1}-\gamma)+ab(\alpha_{1}-\gamma)T_{1}+(D-a)(D-b)u_{R}.
\]
The computed particular solution that corresponds to $(D-a)(D-b)u=ab\alpha_{0}/2-(a+b)(\alpha_{1}-\gamma)$
is $\left(ab\alpha_{0}/2-(a+b)(\alpha_{1}-\gamma)\right)U_{1}$. The
particular solution of 
\begin{equation}
(D-a)(D-b)u=ab(\alpha_{1}-\gamma)T_{1}\label{eq:DaDbT1}
\end{equation}
 is obtained by solving 
\begin{align*}
(D-a)v & =ab(\alpha_{1}-\gamma)T_{1}=(D-a)\left(-b(\alpha_{1}-\gamma)T_{1}\right)+b(\alpha_{1}-\gamma)\\
(D-b)u & =v.
\end{align*}
Because of the way the right hand side of the $(D-a)v$ equation is
rewritten, the particular solution of (\ref{eq:DaDbT1}) may be taken
to be computed as the particular solution of 
\[
(D-b)u=-b(\alpha_{1}-\gamma)T_{1}+b(\alpha_{1}-\gamma)U_{1}'=(D-b)(\alpha_{1}-\gamma)T_{1}-(\alpha_{1}-\gamma)+b(\alpha_{1}-\gamma)U_{1}'.
\]
From the form of the right hand side, we infer that the particular
solution of (\ref{eq:DaDbT1}) is computed to be 
\[
(\alpha_{1}-\gamma)T_{1}-(\alpha_{1}-\gamma)U_{2}+b(\alpha_{1}-\gamma)U_{1}.
\]
Because of the way $f$ was split,
\begin{align}
u^{p} & =\left(ab\alpha_{0}/2-(a+b)(\alpha_{1}-\gamma)\right)U_{1}+(\alpha_{1}-\gamma)T_{1}-(\alpha_{1}-\gamma)U_{2}+b(\alpha_{1}-\gamma)U_{1}+u_{R}\nonumber \\
 & =a\left(b\alpha_{0}/2-\alpha_{1}+\gamma\right)U_{1}-(\alpha_{1}-\gamma)U_{2}+(\alpha_{1}-\gamma)T_{1}+u_{R}\label{eq:computedDaDbup}
\end{align}
is the particular solution of $(D-a)(D-b)u=f$ computed by the factored
form of spectral integration.

The homogeneous solutions computed by the factored form of spectral
integration are 
\begin{align}
\bar{u}^{1} & =\frac{aU_{1}}{2}+\frac{U_{2}}{2}\nonumber \\
\bar{u}^{2} & =\frac{1}{2}+\frac{bU_{2}}{2}.\label{eq:computedDaDbuh}
\end{align}
By observing (\ref{eq:computedDaDbup}) and (\ref{eq:computedDaDbuh}),
we find that 
\[
u_{E}=u^{p}-2\left(b\alpha_{0}/2-\alpha_{1}+\gamma\right)\bar{u}^{1}+\alpha_{0}\bar{u}^{2}.
\]
We may argue as before that even though $u^{p},\bar{\: u}^{1},\bar{\: u}^{2}$
are inaccurate, the factored form of spectral integration solves $(D-a)(D-b)u=f$
with boundary conditions $u(\pm1)=0$ accurately.

\subsection{Condition numbers and accuracy}

\begin{table}
\begin{centering}
\begin{tabular}{|c|c|c|c|}
\hline 
$M$ & error & cond & Bauer\tabularnewline
\hline 
\hline 
$16$ & 5.5e-16 & 3.8e2 & 1.9e1\tabularnewline
\hline 
$32$ & 1.6e-15 & 5.3e3 & 7.3e1\tabularnewline
\hline 
$128$ & 2.9e-15 & 1.2e6 & 1.1e3\tabularnewline
\hline 
$1024$ & 1.1e-13 & 4.8e9 & 7.0e4\tabularnewline
\hline 
$4096$ & 2.5e-13 & 1.5e11 & 2.5e5\tabularnewline
\hline 
\end{tabular}
\par\end{centering}

\caption{Table of errors and condition numbers in the solution of $\left(D^{2}-a^{2}\right)u=f$
with $a=10^{6}$ and $f=-(\pi^{2}+a^{2})\sin\pi y$. The error is
the infinite norm error in the computed $u$. The last two columns
give the standard condition number of the spectral integration matrix
and Bauer's spectral radius.\label{tab:sec3.2-condnums-accuracy-2ndorder}}

\end{table}
Table \ref{tab:sec3.2-condnums-accuracy-2ndorder} shows the errors
in the solution of the linear system $\left(D^{2}-a^{2}\right)u=f$
with $a=10^{6}$ and $f=-(\pi^{2}+a^{2})\sin\pi y$. The errors are
of the order of machine precision when $M=16$ or $M=32$ and grow
only very slowly as $M$ is increased. The version of spectral integration
employed here was that of Section 2.2. However, the results are similar
for the versions in Sections 2.3 or 2.4. 

Table \ref{tab:sec3.2-condnums-accuracy-2ndorder} also shows that
the $2$-norm condition number of the spectral integration matrix
is increasing rapidly. It does converge to a limit as $M\rightarrow\infty$
\cite{CoutsiasHT1996,Greengard1991,Rokhlin1983}, but the limit is
approximately $a^{2}=10^{12}$ (see the last columns of Tables 2 and
3 of \cite{CoutsiasHT1996} for another similar example). The $2$-norm
condition number here has nothing to do with the accuracy of the computed
answer and the fact that it converges to a limit as $M\rightarrow\infty$
is of no consequence.

A more pertinent quantity is Bauer's spectral radius. It is known
that 
\[
\min\kappa_{\infty}\left(D_{1}AD_{2}\right)=\rho\left(|A|\,|A^{-1}|\right)
\]
where the minimum is taken over all non-singular diagonal matrices
$D_{1}$ and $D_{2}$, and $\rho(\cdot)$ is the spectral radius \cite{Bauer1963}\cite[p. 127]{Higham2002}.
Bauer's spectral radius accounts for both row and column scaling.
From table \ref{tab:sec3.2-condnums-accuracy-2ndorder}, this quantity
seems to converge approximately to $a$ and not $a^{2}$ in the limit
$M\rightarrow\infty$. The Green's function corresponding to $\left(D^{2}-a^{2}\right)u=f$
has a scale proportional to $1/a$ (see \cite{ViswanathTobasco2013}).
Therefore, even with row and column scaling we cannot expect a better
condition number than $1/a$.

Although more pertinent, Bauer's spectral radius too fails to explain
the accuracy of computed solution for large $M$ in Table \ref{tab:sec3.2-condnums-accuracy-2ndorder}.
The explanation appears to be that because the spectral integration
matrix is banded and the Chebyshev coefficients of $u=\sin\pi y$
decay rapidly, it is as if only a section of the matrix corresponding
to the lower coefficients is really active. Correspondingly, it may
be noted that the singular vectors corresponding to the largest singular
values are strongly localized within the lowest Chebyshev coefficients.

\begin{table}
\begin{centering}
\begin{tabular}{|c|c|c|c|c|}
\hline 
a & b & M & error1 & error2 \tabularnewline
\hline 
1e+06 & 2e+06 & 1024 & 0.863351 & 0.863351 \tabularnewline
\hline 
1e+06 & 2e+06 & 8192 & 2.14342e-07 & 2.14697e-07 \tabularnewline
\hline 
1e+06 & 2e+06 & 16384 & 1.11927e-09 & 8.68444e-10 \tabularnewline
\hline 
1e+06 & 2e+06 & 131072 & 2.62727e-08 & 3.47769e-08 \tabularnewline
\hline 
\end{tabular}
\par\end{centering}

\caption{Infinite norm errors in the solution of $(D^{2}-a^{2})(D^{2}-b^{2})u=a^{2}b^{2}$
with $u(\pm1)=u'(\pm1)=0.$ The two errors correspond to spectral
integration using the factorizations $(D-a)(D+a)(D-b)(D+b)$ and $(D^{2}-a^{2})(D^{2}-b^{2})$,
respectively. \label{tab:sec3.2-condnum-accuracy-4thorder}}

\end{table}

The situation in Table \ref{tab:sec3.2-condnums-accuracy-2ndorder}
is one extreme. The other extreme is shown in Table \ref{tab:sec3.2-condnum-accuracy-4thorder}.
The solution of the fourth order problem in the latter table develops
boundary layers of size $10^{-6}$. The $2$-norm condition number
for the linear systems is $10^{24}$ and is again totally irrelevant
to the observed accuracy. In this latter table, we never get accuracy
close to machine precision. The observed accuracy implies a loss of
at least $6$ digits. Because the solution develops boundary layers,
the assumption that only the lowest few Chebyshev modes are active
is no longer valid. The solution is of poor quality for $M=1024$
because the grid fails to resolve the boundary layers.

The situation in turbulence simulations is probably in between the
two scenarios. For reasons discussed in the introduction, turbulent
solutions will not develop boundary layers or internal layers as thin
as $\mathcal{O}\left(1/a\right)$. Thus we may summarize the discussion
by saying that the unscaled condition numbers are of no relevance,
that Bauer's spectral radius is more pertinent, and even that quantity
may be unduly pessimistic.

For an illustration of the points made so far in this section, we
turn to a Matlab boundary value solver that uses an integral formulation
\cite{Driscoll2010}. This Matlab implementation does not use Chebyshev
series but works exclusively in the physical domain using quadrature
rules to discretize integral operators. While any of the spectral
integration methods of Section 2 applied to $\left(D^{2}-10^{12}\right)u=-\left(\pi^{2}+10^{12}\right)\sin\pi y,\: u(\pm1)=0$
can find the solution $u=\sin\pi y$ with machine precision using
$M=32$, the physical space Matlab implementation fails to do so.
Here we see one advantage of working using Chebyshev coefficients
instead of in physical space. The Matlab implementation can handle
the problem $\left(D^{2}-a^{2}\right)\left(D^{2}-b^{2}\right)u=f$,
if $a$ and $b$ are both $\mathcal{O}(1)$ which is the simplest
scenario. As a point of comparison, we mention that while Matlab took
$1.5523$ billion cycles for a single solve of that system on a single
core of $2.67$ GHz Intel Xeon 5650, the C/C++ implementation of the
method of Section 2.4 can do the same in $108,600$ cycles. Thus the
C/C++ speed-up is $15,000$. Anecdotal evidence suggests that few
regular Matlab users appreciate the stiff performance penalties, which
can get considerably worse when the hardware configuration includes
multiple cores, high speed network, and accelerators.

\subsection{Accuracy of derivatives}

So far our discussion of accuracy has dealt with the solution. Often
one computes the solution as well as its derivatives. The method of
Section 2.3 expands the solution derivative in a Chebyshev series
while the method of Section 2.2 expands the solution itself in a Chebyshev
series. It may seem that the method of Section 2.3 would yield more
accurate derivatives as it does not entail explicit differentiation.

Muite \cite{Muite2010} has show that to be false for certain examples
provided the derivative is computed carefully (see in particular Figures
7 and 10 of his paper and the associated discussion). Here we introduce
an analogy to strengthen Muite's discussion.

Consider the boundary value problem $\left(D^{2}-a^{2}\right)u=f,\: u(\pm1)=0$
and two forms of spectral integration. The first form expands $u$
in a Chebyshev series (this is the method of Section 2.2) and the
second form expands $D^{2}u$ (this is the method of Section 2.3).
Muite \cite{Muite2010} notes that derivatives may be obtained more
accurately from the first form provided the Chebyshev coefficients
of $u$ are used \emph{directly }to calculate the coefficients of
the derivatives in a suitable trigonometric expansion. The accuracy
will be lost if there is any passage to the physical domain.

An analogy is perhaps useful to clarify this subtle point. Consider
$\left(D^{2}-a^{2}\right)u=f$ with the periodic boundary condition
$u(y+2\pi)=u(y)$. Suppose that $f=\sum_{j=-M}^{M}\hat{f}_{j}\exp(\sqrt{-1}jy)$
is the truncated Fourier expansion of $f$ with $|f_{0}|\approx1$
and rapid decay of Fourier coefficients until $|f_{j}|\approx10^{-16}$
for $j=\pm M$. Thus $f$ is assumed to be represented with accuracy
comparable to machine precision. To clarify the discussion, it is
assumed in addition that $M=100$. 

The Fourier coefficients of $u$ are calculated as $\hat{u}_{j}=-\hat{f}_{j}/\left(a^{2}+j^{2}\right)$.
If $a=1$, the Fourier coefficients $\hat{u}_{j}$ decay from approximately
$1$ for $j=0$ to approximately $10^{-20}$ for $j=\pm M$. All these
Fourier coefficients will have significant digits. Therefore if we
form the Fourier coefficients of $D^{2}u$ in the Fourier domain as
$-j^{2}\hat{u}_{j}$ we will get accurate coefficients from $j=0$
till $j=\pm M$. Furthermore, the computed derivative will be accurate
to machine precision.

Now suppose that $u$ is transformed into the physical domain and
than back to the Fourier domain before differentiation. Although this
is mathematically an identity operation, the rounding errors will
imply that all Fourier coefficients of $u$ smaller than $10^{-16}$
are lost. Therefore $D^{2}u$ computed in this manner will be less
accurate.

Another point of interest emerges if we assume $a=10^{3}$. In this
case, the Fourier coefficients of $u$ decay from approximately $10^{-6}$
for $j=0$ to about $10^{-22}$ for $j=\pm M$. Correspondingly, the
Fourier coefficients of $D^{2}u$ decay from about $10^{-6}$ to about
$10^{-18}$. In the physical domain $D^{2}u$ will have $12$ digits
of relative accuracy which is short of machine precision. Thus the
efficacy of differentiating in spectral domain appears to decrease
as $a$ increases.

\section{Numerical examples}

In this section, we give numerical examples that illustrate the properties
of spectral integration. The examples build up to a discussion of
the properties of the Kleiser-Schumann algorithm and motivate the
form of the algorithm derived in \cite{Viswanath2014}.

\subsection{An example with a boundary layer}

\begin{figure}
\centering{}\includegraphics[scale=0.4]{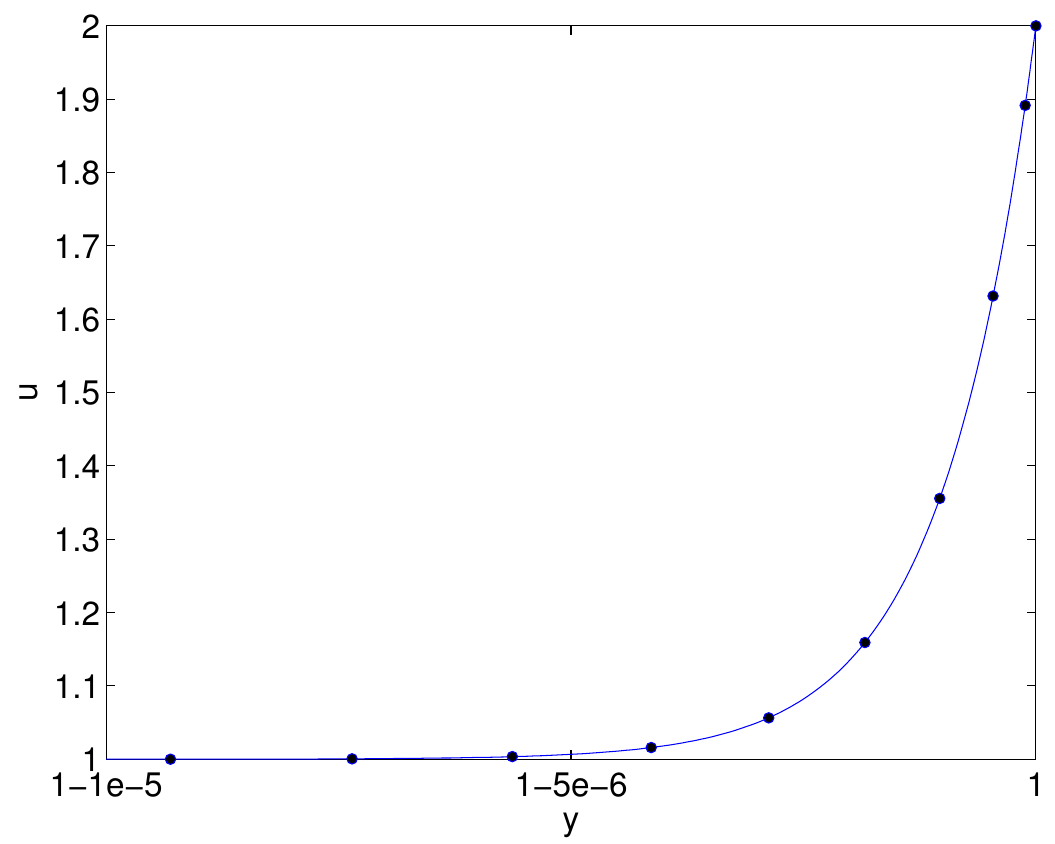}\caption{Solution of $(D^{2}-aD)u=0$ with $a=10^{6}$ and $u(1)=2$ and $u(-1)=1$.
The figure shows the very thin boundary layer near $y=1$. The computed
solution (solid line) is in excellent agreement with the exact solution
(filled markers).\label{fig:2ndorder boundary layer}}
\end{figure}

\begin{table}
\begin{centering}
\begin{tabular}{|c|c|c|c|c|c|}
\hline 
M1 & M2 & M3 & node2 & node3 & error\tabularnewline
\hline 
16 & 1024 & 32 & 0.5 & 0.99999 & 5.80845e-06 \tabularnewline
\hline 
16 & 4096 & 32 & 0.5 & 0.99999 & 4.07361e-11 \tabularnewline
\hline 
32 & 128 & 32 & 0.999 & 0.99999 & 4.49718e-11 \tabularnewline
\hline 
32 & 64 & 32 & 0.9999 & 0.99999 & 4.33247e-11 \tabularnewline
\hline 
32 & 32 & 32 & 0.99995 & 0.99999 & 4.66069e-11 \tabularnewline
\hline 
\end{tabular}
\par\end{centering}

\caption{Solution of $D^{2}-aD=0$ with $u(-1)=1$, $u(1)=2$, and $a=10^{6}$
using a grid with three intervals, which are discretized using $M1+1$,
$M2+1$, and $M3+1$ Chebyshev points, respectively. Nodes 1 and 4
are located at $\pm1$. Node 2 is outside the boundary layer in the
first two rows. The error is in the infinity norm. \label{tab:Solution-of-D2maD}}
\end{table}

Table \ref{tab:Solution-of-D2maD} summarizes spectral integration
of piecewise Chebyshev grids applied to solve $(D^{2}-aD)u=0$ with
$u(-1)=1$, $u(1)=2$, and $a=10^{6}$. This problem develops a boundary
layer at $y=1$; see Figure \ref{fig:2ndorder boundary layer}. It
is evident from the table, that the intervals must be chosen carefully.
The table shows that attempts to get an accurate solution with fewer
than a thousand grid points and just a single interval properly contained
in the boundary layer did not work. The last row of Table \ref{tab:Solution-of-D2maD}
reports a solution with $M_{1}=M_{2}=M_{3}=32$ and an error of $4.7\times10^{-11}$.
In that computation, two intervals are contained inside the boundary
layer. If a single Chebyshev grid is used, $M=8192$ is needed to
get more than ten digits of accuracy.

The example of Table \ref{tab:Solution-of-D2maD} coincides with Example
3 of \cite{GreengardRokhlin1991}. Table 7 of \cite{GreengardRokhlin1991}
reports an error of $2.33\times10^{-11}$ using $20$ intervals and
$M=16$ for each interval (the total number of grid points is $321$).

\subsection{An example with an internal layer}

\begin{figure}
\begin{centering}
\includegraphics[scale=0.3]{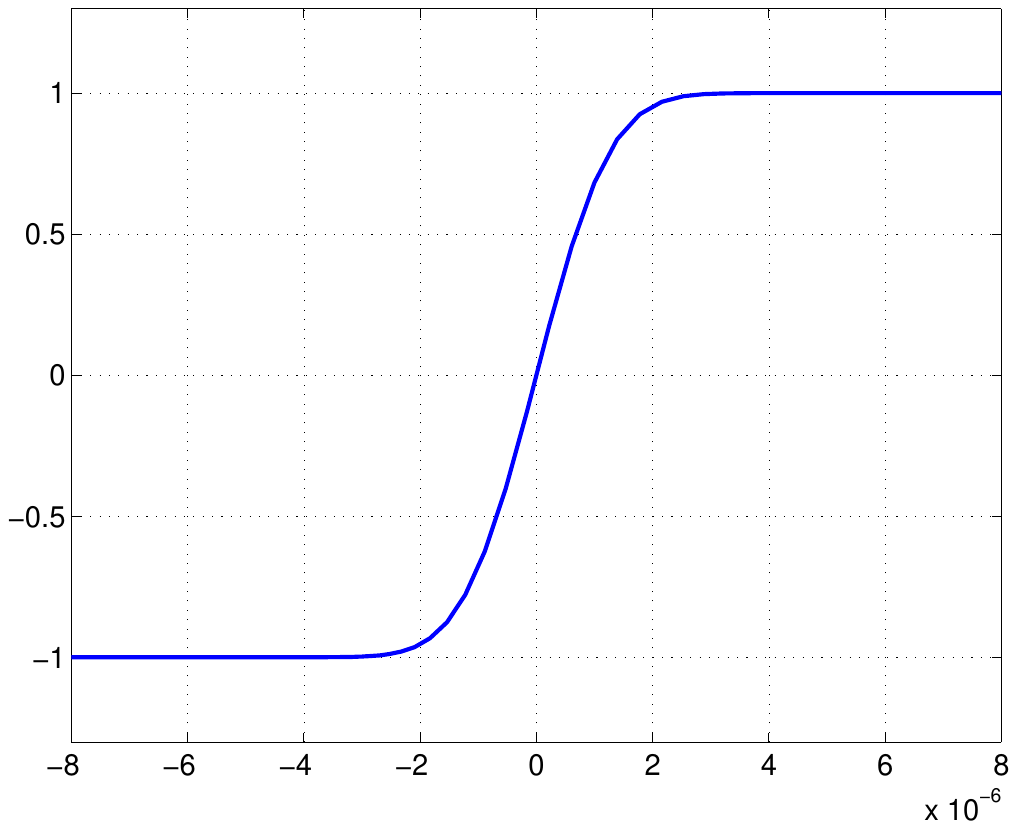}
\par\end{centering}

\vspace{-2cm}

\centering{}\caption{Solution of $u''+a\, y\, u'=0$ with $u(-1)=-1$, $u(1)=1$, and $a=10^{12}$
using differentiation matrices and a piecewise Chebyshev grid. The
plot shows the transition region of the solution.\label{fig:transition-spy-layer}}
\end{figure}

The second example we consider is $\epsilon u''+yu'=0$ with boundary
conditions $u(\pm1)=\pm1$ and $\epsilon=10^{-12}$. The exact solution
of this boundary value problem is given by 
\[
u(y)=-1+\frac{2\int_{-1/\sqrt{2\epsilon}}^{y/\sqrt{2\epsilon}}{\rm e}^{-t^{2}}\, dt}{\int_{-1/\sqrt{2\epsilon}}^{1/\sqrt{2\epsilon}}{\rm e}^{-t^{2}}\, dt}.
\]
The solution has an internal layer at $y=0$ of width approximately
$\epsilon^{-1/2}$ or $10^{-6}$. In Figure \ref{fig:transition-spy-layer},
we show the spy plot of a matrix corresponding to division of $[-1,1]$
into five sub-intervals as well as the transition region of the solution.

\begin{table}
\centering{}%
\begin{tabular}{|c|c|c|}
\hline 
m & node4 & overshoot\tabularnewline
\hline 
32 & $5\sqrt{\epsilon}$ & 3.7e-15\tabularnewline
\hline 
32 & $3\sqrt{\epsilon}$ & 1.2e-08\tabularnewline
\hline 
32 & $7\sqrt{\epsilon}$ & 8.6e-09\tabularnewline
\hline 
24 & $5\sqrt{\epsilon}$ & 1.8e-08\tabularnewline
\hline 
\end{tabular}\caption{Solution of $\epsilon\, u''+y\, u'=0$ with $u(-1)=-1$, $u(1)=1$,
and $\epsilon=10^{-12}$ is shown in Table \ref{fig:transition-spy-layer}.
This table gives the overshoot beyond $[-1,1]$ of the solution computed
using 6 nodes and 5 intervals. Nodes 1, 2, 3, 5, and 6 are fixed at
$-1$, $-8\epsilon^{1/2}$, $-3\epsilon^{1/2}$, $8\epsilon^{1/2}$,
and $1$, respectively. The number of Chebyshev points in each interval
is $m+1$. \label{tab:Solution-of-internal-layer}}
\end{table}

This second example occurs near the end of \cite{CoutsiasHT1996},
where it is reported that mapped Chebyshev points with $M=1024$ compute
the solution with an overshoot of $3\times10^{-4}$. From Table \ref{tab:Solution-of-internal-layer},
we see that the overshoot is reduced to the order of machine precision
using only $161$ grid points. The overshoot is seen to be highly
sensitive to the location of the nodes. The solution plotted in Figure
\ref{fig:transition-spy-layer}b corresponds to the top row of the
table. The solution appears to have around $10$ digits of accuracy.

\subsection{Differentiation errors and the Kleiser-Schumann algorithm}

The two examples above are extreme examples of the efficacy of a piecewise
Chebyshev grid. However, even in the application to the Navier-Stokes
equation, a piecewise Chebyshev grid can give superior resolution
with fewer grid points. Unfortunately, using a piecewise Chebyshev
grid raises the CFL (Courant-Friedrichs-Lewy) number significantly
and imposes an excessive constraint on the time step. Therefore, it
does not prove useful. However, there are situations where the Navier-Stokes
equations are solved for steady solutions and traveling waves without
time integration \cite{Waleffe2003}. Piecewise Chebyshev grids may
prove useful in such applications.

The solution of the Navier-Stokes equation for channel flow or plane
Couette flow involves two derivatives in the wall-normal direction.
The first derivative arises when computing vorticities or the nonlinear
term. The second order derivative arises when computing the right
hand side of the pressure Poisson equation \cite{Viswanath2014,ViswanathTobasco2013}.
Thus the rounding errors that arise during the differentiation operations
may degrade the accuracy of the computed solutions.

\begin{figure}
\begin{centering}
\includegraphics[scale=0.35]{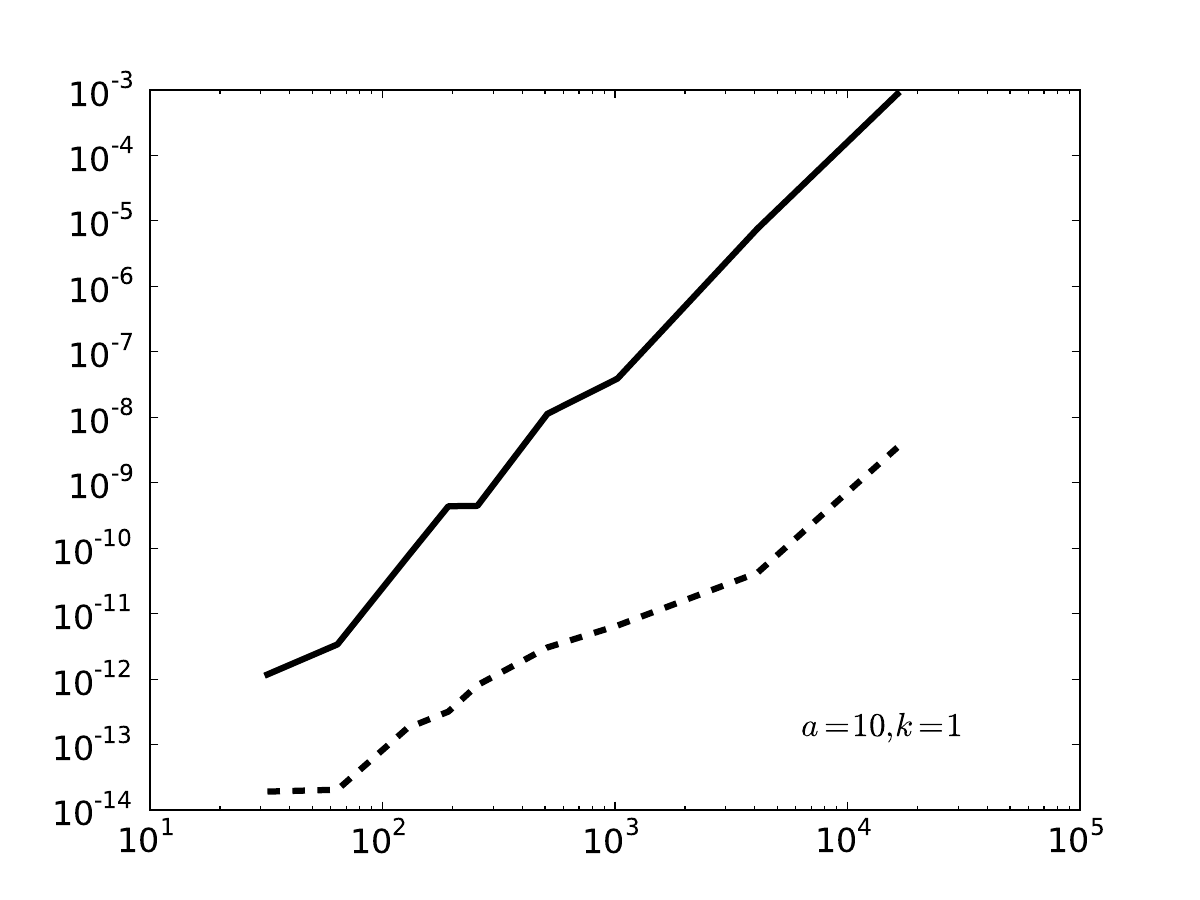}\includegraphics[scale=0.35]{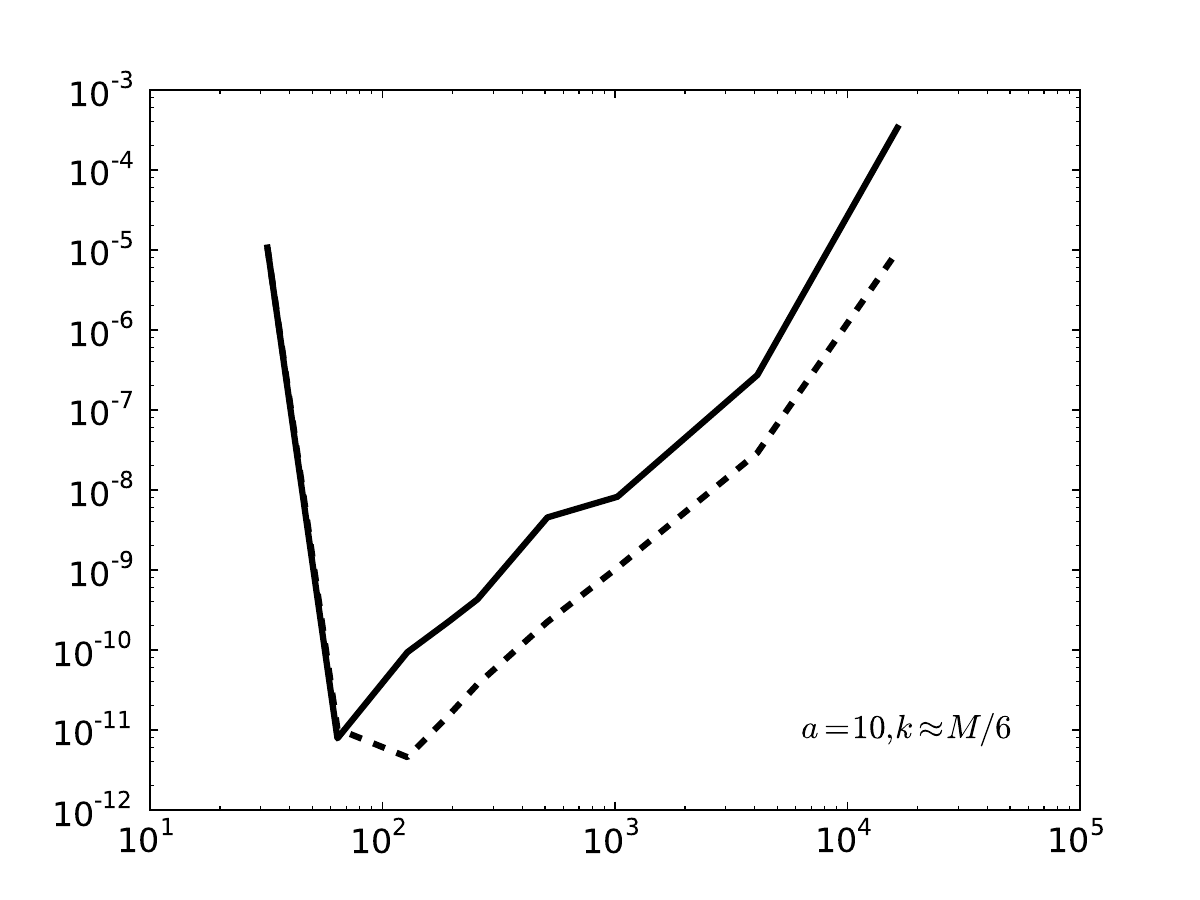}
\par\end{centering}

\begin{centering}
\includegraphics[scale=0.35]{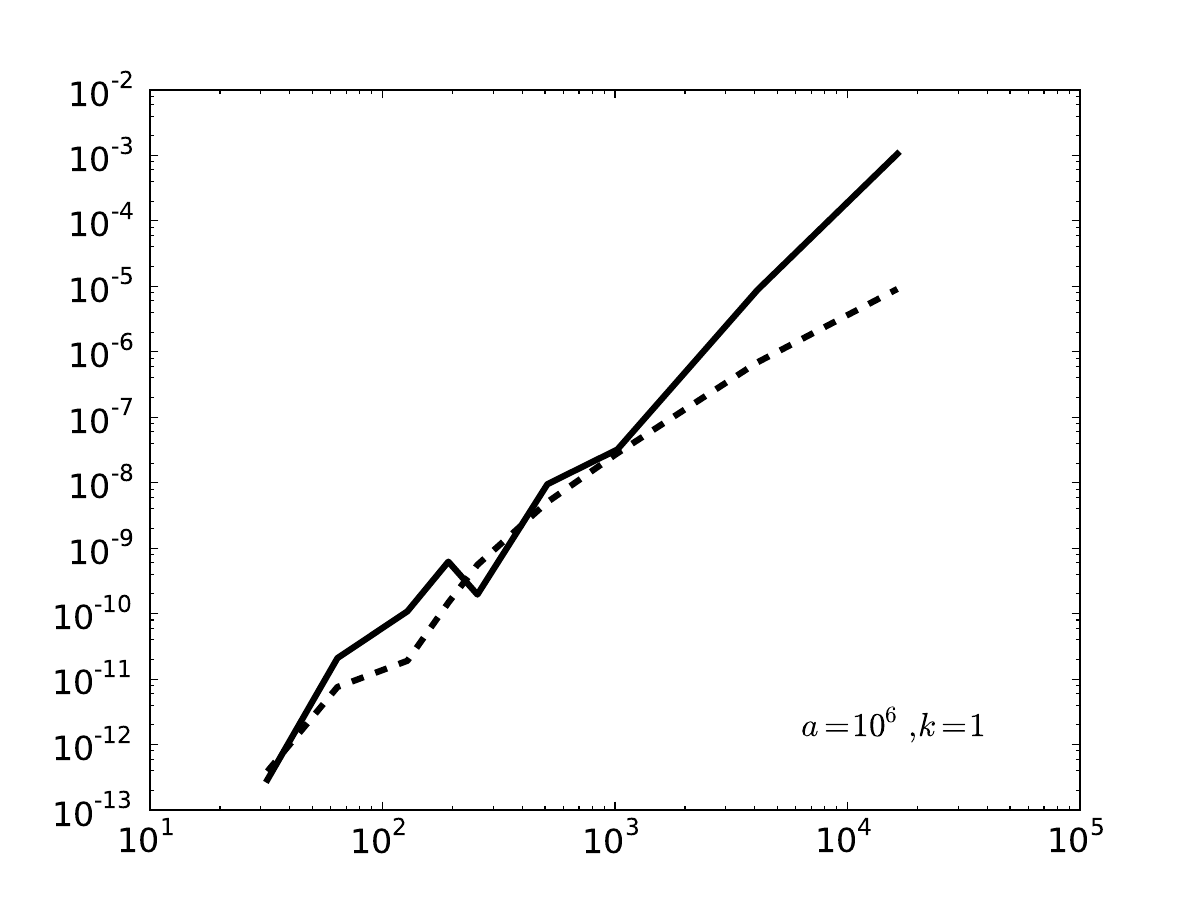}\includegraphics[scale=0.35]{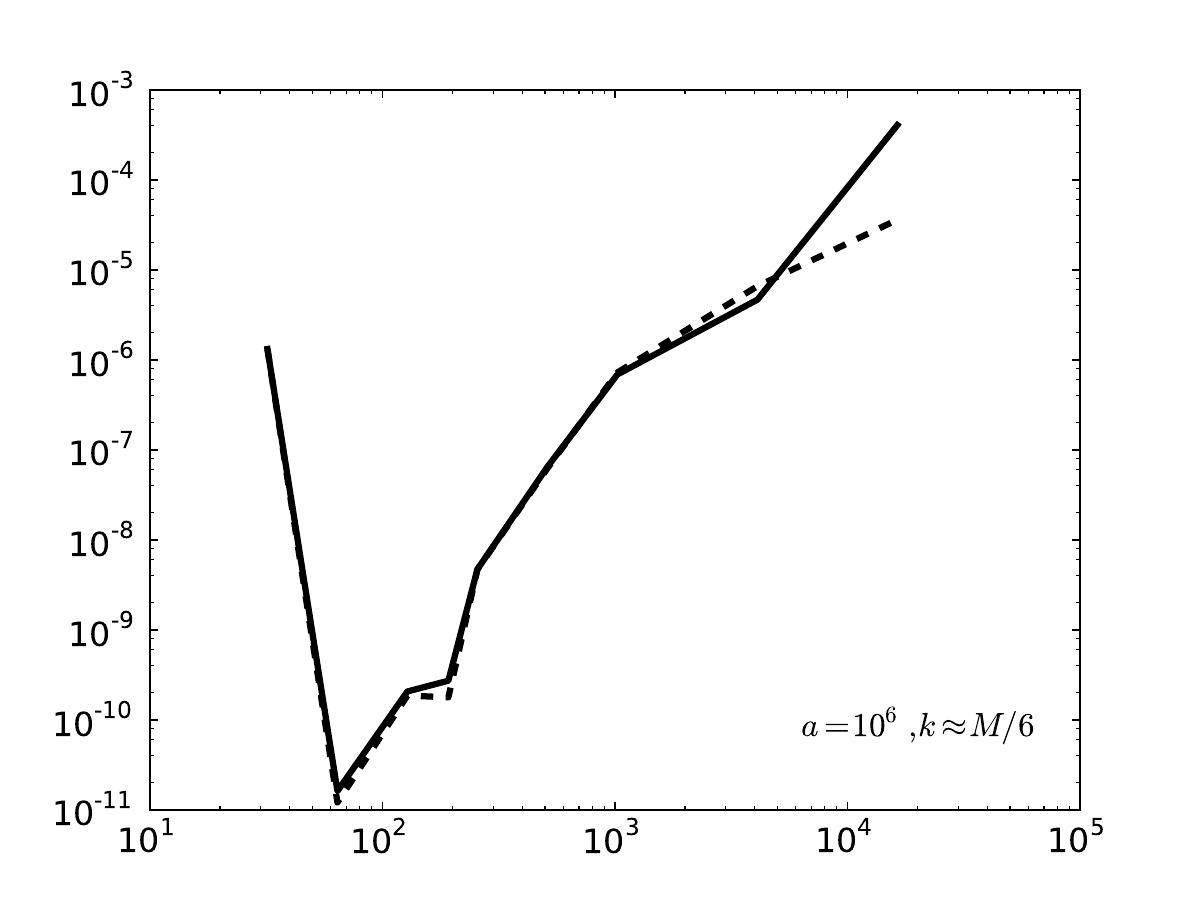}
\par\end{centering}

\caption{Plots of the infinite norm error in the derivative of $du/dy$ of
the solution of $\left(D^{2}-a^{2}\right)u=f$ with boundary condition
$u(\pm1)=0$ against $M$, the number of Chebyshev modes. In the plots
on the left, $f$ is chosen so that the exact solution is $u=\sin\pi y$.
In the plots on the right, $f$ is chosen so that the exact solution
is $f=\sin k\pi y$ with $k\approx M/6$. The dashed line corresponds
to differentiation in the spectral domain without passing to the physical
domain. The solid line corresponds to differentiation  after passing
to the physical domain.\label{fig:secn4.3-diff-errors}}
\end{figure}
 Figure \ref{fig:secn4.3-diff-errors} shows differentiation errors
in four different situations. In each case, the solution of $\left(D^{2}-a^{2}\right)u=f$
is computed using the method of Section 2.1 or 2.2. The difference
between the solid line and the dashed line in the figure is minor.
The derivative corresponding to the solid line is obtained by transforming
the computed $u$ from Chebyshev space to physical space and then
back to Chebyshev space for differentiation using the discrete cosine
and sine transforms. For the dashed line, the computed $u$ is not
transformed to physical space but is differentiated in Chebyshev space
without going to physical space. In both cases, errors in the derivative
are measured in the physical space.

The $a=10$ and $k=1$ plot corresponds best to the calculations in
\cite{Muite2010}. Here we see that the derivative obtained without
passing to physical space is much more accurate. Most of the error
is near the edges and the error would be lower in the energy norm.
If $k\approx M/6$, the advantage of not passing to physical space
is not as great.

The case $a=10^{6}$ is more relevant to turbulence simulations at
high Reynolds number \cite{Viswanath2014}. Here the advantage of
not passing to physical space all but disappears, which is in agreement
with the discussion in Section 3.3. Thus not passing to physical space
appears to imply little advantage in turbulence simulation at high
Reynolds number. In contrast, there is a significant gain in accuracy
if numerical differentiation is avoided entirely by expanding $du/dy$
in a Chebyshev series as shown in \cite{Viswanath2014}.

The heart of the Kleiser-Schumann algorithm \cite{KleiserSchumann1980}
for solving the Navier-Stokes equation in rectangular geometry is
the following set of equations: 
\begin{align*}
\frac{\partial u}{\partial t}+H_{1} & =-\frac{il}{\Lambda_{x}}p+\frac{1}{Re}\left(D^{2}-\alpha^{2}\right)u\\
\frac{\partial v}{\partial t}+H_{2} & =-\frac{\partial p}{\partial y}+\frac{1}{Re}\left(D^{2}-\alpha^{2}\right)v\\
\frac{\partial w}{\partial t}+H_{3} & =-\frac{in}{\Lambda_{z}}p+\frac{1}{Re}\left(D^{2}-\alpha^{2}\right)w.
\end{align*}
Here $2\pi\Lambda_{x}$ and $2\pi\Lambda_{z}$ are the extent of the
domain in the wall parallel directions, and $u(y),v(y),w(y)$ are
the coefficients of $\exp(il/\Lambda_{x}+in/\Lambda_{z})$ in the
Fourier expansion of the velocity field. In addition, $\alpha^{2}=l^{2}/\Lambda_{x}^{2}+n^{2}/\Lambda_{z}^{2}$
and $H_{1},H_{2},H_{3}$ are nonlinear advection terms. The incompressibility
condition is $(il/\Lambda_{x})u+\partial v/\partial y+(in/\Lambda_{z})w=0$.
Using the incompressibility condition, we get an equation for pressure:
\[
\left(D^{2}-\alpha^{2}\right)p=-\frac{il}{\Lambda_{x}}H_{1}-\frac{\partial H_{2}}{\partial y}-\frac{in}{\Lambda_{z}}H_{3}.
\]
After time discretization, the Kleiser-Schumann algorithm solves these
equations for the Chebyshev series of $u,v,w,p$. The type of errors
discussed here occur when $H_{1},H_{2},H_{3}$ are computed for use
in the next time step. In this computation, $u$ and $w$ are differentiated
with respect to $y$. The nonlinear term $H_{2}$ is once again differentiated
with respect to $y$ when solving for pressure. 

The derivatives of $u,w$ with respect to $y$ may be formed without
passing to the physical domain. Such an implementation would be more
accurate at low Reynolds numbers, for which the parameter $a$ in
$\left(D^{2}-a^{2}\right)u=f$ is small. But the advantage of not
passing to the physical domain diminishes as the Reynolds number increases
as evident from Figure \ref{fig:secn4.3-diff-errors}.

\section{Conclusions}

In this article, we have derived many different versions of spectral
integration for solving linear boundary value problems. The treatment
of boundary conditions is uncoupled from finding a particular solution
in each of these versions. As a result, the various forms of spectral
integration can be combined in many ways. 

The accuracy of spectral integration has a number of subtle features.
In Section 3, we discussed the reason one can get away with not having
to resolve boundary layers present in the Green's function but not
the solution as well as condition numbers and the accuracy of derivatives.
The discussion is illustrated in Section 4 using numerical examples
and motivates the new versions of Kleiser-Schumann and Kim-Moin-Moser
algorithms derived in \cite{Viswanath2014}. The new version of Kleiser-Schumann
has been used to carry out a computation with $10^{9}$ grid points
at $Re_{\tau}=2380$, which appears to be the highest Reynolds number
reached in fully resolved simulations of wall bounded turbulence.

\section{Acknowledgements}

Thanks to Fabian Waleffe for comments and suggestions. Thanks as well
to Hans Johnston and Benson Muite. This research was partially supported
by NSF grants DMS-1115277 and SCREMS-1026317.

\bibliographystyle{plain}
\bibliography{references}

\end{document}